\documentclass[12pt]{amsart}
\usepackage{amscd,amsthm,amsmath} 
\newtheorem{thm}{Theorem}[section]
\newtheorem{cor}[thm]{Corollary}
\newtheorem{lem}[thm]{Lemma}
\newtheorem{prop}[thm]{Proposition}
\newtheorem{examp}[thm]{Example}
\newtheorem{ques}[thm]{Question}
\newtheorem{conj}[thm]{Conjecture}

\theoremstyle{definition}
\newtheorem{defn}[thm]{Definition}
\theoremstyle{remark}
\newtheorem{rem}[thm]{Remark}
\numberwithin{equation}{section}

\begin{document}

\title{Phantom elements and its Applications}%
\author{Jianzhong Pan }%
\address{Institute of Math.,Academia Sinica ,Beijing 100080, China }%
\email{pjz@math03.math.ac.cn}%
\author{Moo Ha Woo}
\address{ Department of Mathematics Education , Korea University ,
 Seoul , Korea}

\thanks{The first author is partially supported by the NSFC project
19701032 and ZD9603 of Chinese Academy of Science and the second
author wishes to acknowledge the financial supports of the Korea
Research Foundation made in the
program year of (1998) and TGRC 99}%

\date{Aug. 30, 2000}%
\begin{abstract}
 In our previous work\cite{panwoo} , a relation between Tsukiyama problem
  about self homotopy equivalence  was found by using a generalization of
phantom map. In this note ,  fundamental result is established for
such a generalization. This is the first time one can deal with
phantom maps to space  not satisfying finite type condition.
Application to Forgetful map is also discussed briefly.
\end{abstract}
\maketitle
 \section{Introduction} \label{S:intro}
The main aim of this paper is to study phantom map , its generalization
and  applications. After the discovery of the first example of phantom
 map by  Adams and Walker\cite{aw} , theory of phantom map receives a
  lot of attention . The main aim of these previous studies is however
  to understand it , e.g., the computation and the properties of phantom
   maps. The first application of theory of phantom maps was given by
    Harper and Roitberg\cite{roitharp},\cite{roit2}  who applied it to
    compute $SNT(X)$ and $Aut(X)$.
     Recently applications are also found by Roitberg\cite{roit}
      and Pan\cite{pan} where several conjectures of McGibbon were settled.
       On the other hand , a remarkable connection was established by
        Pan and Woo\cite{panwoo}  between  Tsukiyama problem about
self homotopy equivalence  and a generalization of phantom map. A
byproduct of this connection is that a special case of Tsukiyama
 problem is almost equivalent to the famous Halperin conjecture in
 rational homotopy theory\cite{dupont}.

A well known characterization of map between nilpotent spaces of
finite type to be phantom map is the following
\begin{thm}\label{T:fund}
Let $X,Y$ be nilpotent CW complexes of finite type with $Y$
1-connected and $f:X \to Y $ be any map. Then the followings are
equivalent
\begin{itemize}
\item{$f$ is a phantom map }
\item{$e\circ f \simeq *$ where $e : Y \to \hat{Y}$ is the profinite
 completion }
\item{$f \circ \tau \simeq *$  where $\tau:X_{\tau} \to X$ is the homotopy fiber
 of the rationalization}
\end{itemize}
\end{thm}

On the other hand , in our previous paper\cite{panwoo}, we
generalized the concept of phantom map to that of phantom element
and announced a theorem characterizing an element to be a phantom
element which generalizes Theorem\ref{T:fund} . In this paper we will generalize further so that we can deal with space which is not of finite type.
\begin{thm}\label{T:fundamental}
Let $X$ be nilpotent CW  complex  of finite type , $Y$ be
1-connected such that $\pi_n(Y)$ is reduced group for $n\geq 2$
and $g:X \to Y $ be any map.Then the followings are equivalent:
\begin{itemize}
\item{$\alpha \in \pi_j(map_*(X,Y);g)$ is a phantom element }
\item{$(e_*)_{\#}(\alpha)=0$ where $(e_*)_{\#} :
 \pi_j(map_*(X,Y);g) \to  \pi_j(map_*(X,\hat{Y});\hat{g})$}
\item{$(\tau^*)_{\#}(\alpha)=0$  where $(\tau^*)_{\#}:\pi_j(map_*(X,Y);g)
\to  \pi_j(map_*(X_{\tau},Y);g_{\tau})$}
\end{itemize}
\end{thm}
Note that the assumption that $Y$ is 1-connected is not a real
restriction  by an observation of Zabrodsky\cite{za}.   We will
give a complete proof of this theorem in this paper.

As an application we have(Corollary of Proposition \ref{T:formal})

\begin{cor}
Let $P$ be 1-connected finite dimensional CW complex or such that
$H^*(P,\mathbb{Z}_p)$ is locally finite over $A_p$ for each prime
$p$  and be of type $F_0$ .  Assume further that $\pi_nBaut(P)$
 is reduced group for  $n\geq 2$. Assume further that $P_{(0)}$
satisfying one of the following.
\begin{itemize}
\item{$P$ is rationally equivalent to K\"ahler manifold}
\item{$H^*(P;Q)$ as an algebra has at most 3 generators}
\item{$P$ is rationally equivalent to $G/U$ where $G$ is a compact Lie group and $U$ is a closed subgroup of maximal rank}
\end{itemize}
Then for all $m \geq 1$ , $H$ and every principal $K(H,2m)$-bundle
with total space homotopy equivalent to $P$ , Forgetful map is injective.
\end{cor}

 The organization of this paper is as follows. In section\ref{S:phan}
Theorem\ref{T:fundamental} will be proved.
  The applications to Forgetful map will be discussed in
  section\ref{S:appl}.
  In this paper,  We will  use the
following notations:
\begin{itemize}
\item{$H$ will denote a finitely generated abelian group}
\item{$map(X,Y)$ is the space of continuous mappings from $X$ to $Y$}
\item{$map_*(X,Y)$ is the  subspace of pointed mappings from $(X,x_0)$
to $(Y,y_0)$}
\item{$l:X \to X_{(0)}$ is the rationalization}
\item{ Let $\tau:X_{\tau} \to X$ be the homotopy fiber of $l$.
Then $X_{\tau} \overset{\tau}{\to} X \to X_{(0)}$ is a cofibration
up to homotopy }
\item{$e_p : Y \to \hat{Y}_{\mathbb{Z}_{p^{\infty}}}$ is Bousfield-Kan's p-completion. Let $\hat{Y}=\underset{p}{\prod} \hat{Y}_{\mathbb{Z}_{p^{\infty}}}$ and $e=(e_2, e_3, \cdots):Y \to \hat{Y}$.
Let $Y_{\rho}$ be the homotopy fiber of $e$ }
\end{itemize}
The readers should refer to \cite{panwoo} for all the other
notations which have not been explained here.

In concluding the Introduction , we 'd like to give the following
\begin{conj}
The condition that  $\pi_nY$ is reduced group in this paper can be removed.
\end{conj}

\section{Phantom elements}\label{S:phan}
Let's begin with  definition.
\begin{defn}
Let spaces $X$ be a CW  complex, $Y$ be a space  and $g:X \to Y$
any map. Then an element $\alpha \in\pi_j(map_*(X,Y);g) $ is
called a $g$-phantom element if $(i_n^*)_{\#}(\alpha)=0$ for all
$n \geq 0$ where $(i_n^*)_{\#}: \pi_j map_*(X,Y)\to \pi_j
map_*(X^n,Y)$ is the homomorphism induced by the inclusion $i_n
:X^n \to X$. Denoted by
\[
Ph_j^g(X,Y)=\{\alpha \in\pi_j(map_*(X,Y);g)| \alpha \text{ is a
$g$-phantom element}  \}
\]
\end{defn}

Obviously if $g=$constant and $j=0$ , then $\alpha$ is a
$g$-phantom element iff it represents the homotopy class of a map
which is a phantom map.

The main aim of this section is to prove
Theorem\ref{T:fundamental}. Before that , let's give some results
necessary to the proof .
\begin{lem}\label{T:lemma}
Let  $Y$ be 1-connected such that $\pi_n(Y)$ is reduced group for
 $n\geq 2$. Then
\begin{itemize}
\item{$\pi_n(\hat{Y})=\underset{p}{\prod} Ext(\mathbb{Z}_{p^{\infty}}, \pi_n(Y))$}
\item{For $W$ a finite CW complex, $e_*:\pi_j map_*(W,Y)_f \to \pi_jmap_*(W,\hat{Y})_{\hat{f}}$ is injective}
\end{itemize}
\end{lem}
\begin{proof}
The first statement follows from the fact that $Hom(\mathbb{Z}_{p^{\infty}}, B)=0$ for a reduced group since otherwise there will be nontrivial divisible subgroup in $B$.

To prove the second statement , note that the induced map $\pi_n(Y_{\rho}) \to \pi_n(Y)$ is trivial since $\pi_n(Y)$ is reduced group and   $\pi_n(Y_{\rho})$ is rational thus divisible by the arithmetic square Theorem\cite{dror}.
It follows that $e_*:\pi_n(Y) \to \pi_n(\hat{Y})$ is injective and an easy induction
argument shows what we want for $j \geq 1$. For $j=0$, it can be proved by an argument similar to that of Theorem 2.5.3 of \cite{hilton}.
\end{proof}
\begin{prop}
Let $X$ be   nilpotent space   and $Y$ be 1-connected such that
$\pi_n(Y)$ is reduced group for  $n\geq 2$. Then the followings
hold:
\begin{itemize}
\item{$[\Sigma^n X_{(0)}, \hat{Y}]=* ,
\tilde{H}^n(X_{(0)},\pi_i(\hat{Y}))=0 $  for all $n,i \geq 0$}
\item{$[\Sigma^n X_{\tau}, Y_{\rho}]=* ,
\tilde{H}^n(X_{\tau},\pi_i(Y_{\rho}))=0$ for all $n,i \geq 0$}
\end{itemize}
\end{prop}
\begin{proof}
That $\tilde{H}^n(X_{(0)},\pi_i(\hat{Y}))=0 $  follows from the fact that
\[
Hom(A,B)=0 ,  Ext(A,B)=0
\] for ratinal group $A$ and $B=\underset{p}{\prod} Ext(\mathbb{Z}_{p^{\infty}}, B')$

Then the equation
\[
[\Sigma^nX_{(0)}, \hat{Y}]=\underset {\leftarrow}{\lim}_n
[\Sigma^nX_{(0)}, \hat{Y}^{(n)}]
\]
implies the first statement.

  The
  equation about cohomology in second statement is  true since
  $\pi_n(Y_{\rho})$ is rational while the proof of another equation is
   similar to that as in the first statement.

\end{proof}
\begin{prop}\label{T:4iso}
Let $X$ be   nilpotent space   and $Y$ be 1-connected such that
 $\pi_n(Y)$ is reduced group for  $n\geq 2$. Then the
followings hold :
\begin{itemize}
\item{$\tau^*: map_*(X,\hat{Y}) \overset{w}{\simeq} map_*(X_{\tau},\hat{Y})$}
\item{$\rho_*:map_*(X_{(0)},Y_{\rho}) \overset{w}{\simeq} map_*(X_{(0)},Y)$ }
\item{$e_*:map_*(X_{\tau}, Y) \overset{w}{\simeq} map_*(X_{\tau}, \hat{Y})$  }
\item{$l^*:map_*(X_{(0)}, Y_{\rho}) \overset{w}{\simeq} map_*(X, Y_{\rho})$  }
\end{itemize}
\end{prop}
\begin{proof}
The first and the last statements follow from the last Proposition
and the
 well known Zabrodsky Lemma.
The second and third statements follow from
 a $\lim^1$ argument for $j \geq 0$
\[
* \to \underset{\leftarrow}{\lim}_n ^1\pi_{j+1}map_*(Z,E_n) \to
\pi_jmap_*(Z,E) \to
\underset{\leftarrow}{\lim}_n\pi_jmap_*(Z,E_n)\to *
\]
where in the second statement, $Z=X_{(0)}$ and $E_n$ is the $n-th$
term in the Postnikov-Moore tower of the map $\rho:Y_{\rho}\to Y$
while in the third statement , $Z=X_{\tau}$ and $E_n$ is the
$n-th$ term in the Postnikov-Moore tower of the map $e:Y \to
\hat{Y}$. In both case the sequence $\{\pi_jmap_*(Z,E_n)\}$ is a
sequence consisting of isomorphisms and thus the $\lim^1$ is
trivial and the wanted isomorphisms follows immediately.
\end{proof}
\begin{prop}\label{T:comp}
Let $X$ be   nilpotent space   and $Y$ 1-connected such that $\pi_n(Y)$
 is reduced group for  $n\geq 2$. Let $g:X
\to \hat{Y}$ be any map. Then
\[
Ph_j^g(X,\hat{Y})=*
\]
\end{prop}
\begin{proof}
$Ph_j^g(X,\hat{Y})$ is the $\lim^1$ of a sequence of compact
groups and continuous homomorphisms which is will known to be
trivial.
\end{proof}
\begin{prop}\label{T:rho}
Let $X,Y$ be two nilpotent spaces  with $Y$ 1-connected such that $\pi_n(Y)$
 is reduced group for  $n\geq 2$. Then the
followings hold :
\begin{itemize}
\item{$ map_*(X_{\tau},Y) \overset{w}{\simeq} map_*(X,\hat{Y})$}
\item{$ map_*(X_{(0)},Y) \overset{w}{\simeq} map_*(X,Y_{\rho})$}
\end{itemize}
\end{prop}
\begin{proof}
This is an easy consequence of the Proposition above.
\end{proof}

\begin{proof}[Proof of Theorem\ref{T:fundamental}]
The equivalence between the last two statements follows directly from
 the following commutative diagram where the bottom horizontal homomorphism
 and the right side vertical homomorphism are  isomorphisms by
  Proposition\ref{T:4iso}   .
\[
\begin{CD}
\pi_j map_*(X,Y) @>(\tau^*)_{\#}>>  \pi_j map_*(X_{\tau},Y) \\
@V(e_*)_{\#}VV @V(e_*)_{\#}VV \\ \pi_j map_*(X,\hat{Y})
@>(\tau^*)_{\#}>> \pi_j map_*(X_{\tau},\hat{Y})
\end{CD}
\]
Now assume the first statement, then we have  $(i_n^*)_{\#}(e_*(\alpha))=0$
 for all $n \geq 0$. It follows from Proposition \ref{T:comp}
 that $e_*(\alpha))=0$.
The proof of another direction is similar to that in \cite{panwoo} using Lemma\ref{T:lemma} instead of Sullivan's origional result which is stated only for space of finite type.
\end{proof}
\begin{rem}
It is easy to see that the above proof follows the same pattern as
that given by Oda and Shitanda\cite{os}. We give a prove here
because Oda informed us that there were gaps in their proof and
he don't know if the result is true or not. The similar proof
applies also to the equivariant case which will be discussed in
future publication.
\end{rem}
As noted in \cite{panwoo} , the natural question related to the
application of phantom element to the forgetful map is
\begin{ques}
For two maps $f,g:X \to Y$,what is the relation between
$Ph_j^g(X,Y)$ and $Ph_j^f(X,Y)$?
\end{ques}

\begin{prop}\label{T:fprop}
Let $X,Y$ be nilpotent CW  complexes such that \[[\Sigma^j
X_{\tau},Y]=[\Sigma^{j+1}X_{\tau},Y]=0\] If $g:X \to Y$ is a
phantom map,then we have
\[
 Ph^g_j(X,Y)=\pi_j(map_*(X,Y);g)
\]
\end{prop}

\begin{proof}
The proof is the same as that in \cite{panwoo} .
\end{proof}

In our application  we have to be able to compute  $Ph^g_j(X,Y)$.
 Before giving this kind of result, recall that a CW complex is
  called unstable if all the attaching maps vanish under suspension.
  It is Baues \cite{baues} who noted the following
 which is dual to Zabrodsky's integral approximation.
\begin{thm}
Let $X$ be 1-connected CW complex. Then there is an unstable complex and a rational equivalence $h: \bar{X} \to X$.
\end{thm}
\begin{rem}
Let $X$ be an unstable CW complex. Then it is easy to prove that
$Ph_j^g(X,Y)=*$ for any map $g:X \to Y$.
\end{rem}

\begin{prop}\label{T:calc}
Let $X$ be a 1-connected CW complex and  $Y$ 1-connected such that
$\pi_n(Y)$
 is reduced group for $n\geq 2$. Suppose further that the component of
 $map_*(X, \hat{Y})$ consisting constant map is
  weakly contractible and $g:X \to Y$ is a phantom map. Then
\[
 Ph^g_j(X,Y)=\pi_j(map_*(X,Y);g)=
 \underset{k>0}{\prod}H^k(X, \pi_{k+j+1}(Y_{\rho}))
\]
\end{prop}
\begin{proof}
As first noted by Oda and Shitanda , similar proof as in that of
Theorem B of \cite{za} leads to the following homotopy fibration
\[
\underset{g}{\bigcup}map_*(X,Y)_g \to map_*(\bar{X},Y)_* \to
map_*(\bar{X},\hat{Y})_*
\]
where the union is over phantom maps $g$. On the other hand ,
different components of $\underset{g}{\bigcup}map_*(X,Y)_g$ are
homotopy equivalent since $\underset{g}{\bigcup}map_*(X,Y)_g$ is
the homotopy fiber of a map between two connected spaces.  It
follows that
\[
 Ph^g_j(X,Y)=\pi_j(map_*(X,Y);const)=\pi_j(map_*(X,Y_{\rho});const)=
 \]
 \[
 =[\Sigma^{j-1}X, \Omega Y_{\rho}]=
 \underset{k>0}{\prod}H^k(X,\pi_{j+k+1}(Y_{\rho}))
 \]
\end{proof}
We are ready to state results related to the applications. Before
that we have another definition
\begin{defn}
Let $A_p$  be the modp Steenrod algebra. An unstable module $M$
over $A_p$ is called locally finite iff , for any $x \in M$, only
finite  elements of $A_p$  can acts nontrivially on $M$.
\end{defn}
\begin{examp}
Let $P$ be a space such that $H^*(P,\mathbb{Z}_p)$ is locally
finite over $A_p$ . Then so is $\Omega P$. In particular, if $P$
is finite CW , then $H^*(\Omega P,\mathbb{Z}_p)$ is locally finite
over $A_p$.
\end{examp}
\begin{thm}\label{T:EMphan}
Let $X=K(H,m+2)$ , $Y=Baut(P)$ such that $\pi_n(Y)$
 is reduced group for  $n\geq 2$ and $g:X \to Y $ is any map where
$P$ is 1-connected finite dimensional CW complex or such that
$H^*(P,\mathbb{Z}_p)$ is locally finite over $A_p$ for each prime
$p$ and $m \geq 1$. Then for $j \geq 1$
\[
 Ph^g_j(X,Y)= \pi_j(map_*(X,Y);g)
 =[\Sigma^j X,Y_{\rho}]
\]
\end{thm}
\begin{proof}
The proof is the same as that of the corresponding result in
\cite{panwoo} using results of Zabrodsky and Miller\cite{miller} or Theorem 8.8 in
\cite{schwartz}.
\end{proof}
Similarly we have
\begin{thm}\label{T:liephan}
 Let $X=BG$, $Y=Baut(P)$ such that
$\pi_n(Y)$
 is reduced group for  $n\geq 2$. and $g:X \to Y$ is a phantom map where
 $G$ is a connected compact Lie group and  $P$ is 1-connected finite
dimensional CW complex or such that $H^*(P,\mathbb{Z}_p)$ is
locally finite over $A_p$ for each prime $p$. Then for $j \geq 1$
we have
\[
Ph_j^g(X,Y)=\pi_j(map_*(X,Y);g)= [\Sigma^jX,Y_\rho]
\]
\end{thm}

\section{Application to the forgetful map}\label{S:appl}

 Given
a principal G-bundle $\pi:P \to B$, Let
\[
aut^G(P)=\{g | g:P \to P  \text{  is a G-equivariant homotopy
equivalence}\}
\]
and
\[
aut(P)=\{g | g:P \to P  \text{  is a homotopy equivalence}\}
\]
There is a natural map $f:aut^G(P)\to aut(P)$. Let
\[
Aut^G(P)= \pi_0(aut^G(P))
\]
  and
\[
  Aut(P)= \pi_0(aut(P))
\]
Then the map $f$ induces   a map
\[
F:Aut^G(P) \to Aut(P)
\]
which is called a Forgetful map by Tsukiyama. The question posed
by Tsukiyama in \cite{dwk} is the following
\begin{ques} \label{T:ques2}
Is the forgetting  map $F$   injective?
\end{ques}
One of the main results in \cite{panwoo} is the following
\begin{thm}\label{T:main5}
Let $\pi:P \to B$ be a principal $G$-bundle. Then there is an
exact
 sequence
\[
\pi_1 aut(P) \overset{\delta}{\to} \pi_1(map_*(BG,Baut(P)),c) \to
Aut^G(P) \overset{F}{\to} Aut(P)
\]
where  $c:BG \to Baut(P)$ is determined by the principal bundle.
\end{thm}

Combined with results in \cite{panwoo} , we have
\begin{prop}\label{T:formal}
Let $P$ be as in Theorem\ref{T:liephan} . If
$$\bigoplus_{i>1}\pi_{2i}(map(P_{(0)},P_{(0)});id)=0$$  then for
all $m \geq 1$ , finitely generated abelian group $H$ and every
principal $K(H,2m)$-bundle with total space homotopy equivalent to
$P$ ,the associated Forgetful  map is injective  .
\end{prop}
We have also similar results for $K(H,2m+1)$ or $G$ bundle where
$G$ is a connected compact Lie group which will be omitted.

 Unlike that in \cite{panwoo}, there are no complete results if
 group \newline $\pi_1(map_*(BG,Baut(P)),c)$ is nontrivial although we
 know that it is still uncountable since group $\pi_1aut(P)$
 itself may be uncountable too. Thus same results as in
 \cite{panwoo} can be obtained if $\pi_1aut(P)$ is countable. This
 is so if $P=\Omega P'$  where $P'$ is finite complex. An
 interesting question is
 \begin{ques}
Study the map $\delta$ in the exact sequence of Theorem
\ref{T:main5}. Is it possible that $Image(\delta)$ is always
countable group.
 \end{ques}
----------------------------------------------------------------

----------------------------------------------------------------

\end{document}